\documentclass[12pt,oneside]{amsproc}
\usepackage[T2A]{fontenc}
\usepackage[cp1251]{inputenc}
\usepackage{amssymb}
\usepackage{amsthm}

\title{О-фреймы для операторов в банаховых пространсвах}
\author{Олег Рейнов}
\address{ St. Petersburg State University,
Saint Petersburg, RUSSIA.}
\email{orein51@mail.ru}

\thanks{${ }$ Key words: 46B28  Approximation of operators; bounded approximation property.}

\begin{document}

                          $$ {} $$
\maketitle


 В работе введено новое понятие О-фрема для оператора в банаховых прост\-ранствах.
Исследованы ограниченно полные и натягивающие О-фреймы.
Показано, что оператор имеет О-фрейм тогда и только тогда,
когда 1)\, он факторизуется через пространство с базисом, и, 
если он действует из сепарабельного пространства,
то тогда и только тогда, когда  2)\, он обладает свойством ограниченной ап\-про\-кси\-мации.
Среди прочего, аналогичный факт установлен для введенных здесь же
БО-фреймов (безусловных операторных фреймов).

 Также, например, установлено, что всякий оператор имеющий 
ограниченно полный и натягивающий О-фрейм, слабо компактен.
 Приведены  и другие применения. Одно из них: дан отрицательный ответ на вопрос
 Марии  Скопиной о существовании в $L_1(0,1)$
 безусловной системы представления.

\vskip 0.3cm

\centerline{\bf \S1. О-фреймы}

 \vskip 0.23cm
Мы используем стандартные обозначения и терминолргию
теории операторов в банаховых пространствах (см., например, [2]).      
\smallskip

{\bf Определение 1.1.}\,
Пусть $T\in L(X,W),$ $(x'_k)_{k=1}^\infty\subset X^*, (w_k)_{k=1}^\infty\subset W.$  
Мы говорим, что $\mathcal F:=((x'_k)_{k=1}^\infty, (w_k)_{k=1}^\infty)$
есть О-фрейм для $T,$ если для каждого $x\in X$
ряд $\sum_{k=1}^\infty \langle  x'_k, x\rangle w_k$ сходится в $W$ и
$$
  Tx= \sum_{k=1}^\infty \langle  x'_k, x\rangle w_k, \ \ x\in X.
$$
Если существует О-фрейм для $T,$ то мы говорим, что $T$ имеет О-фрейм.
\smallskip

{\bf Примеры,}\
1)\,
Пусть $\Delta: l_\infty\to l_1$ --- диагональный оператор с диагональю $(\delta_k)\in l_1.$
Тогда $\Delta x=\sum \langle  e_k, x\rangle e_k,$ где $(e_k)$ --- последовательность единичных векторов в $l_p$
(в нашем случае, в $l_1$), $x=(x_k)\in l_\infty.$

2)\,
Если в $X$ есть базис Шаудера $(f_k)_{k=1}^\infty,$
то для любого $x\in X$ имеем: $x=\sum_{k=1}^\infty \langle  f'_k, x\rangle f_k,$ где
$(f'_k)$ --- биортогональная к $(f_k)$ система.
Следовательно, для любого банахова пространства $W$ и каждого оператора
$T:X\to W$
$$
 Tx= \sum_{k=1}^\infty \langle   f'_k, x\rangle Tf_k,\ \ x\in X.
$$

3)\, 
Если $W$ имеет базис Шаудера, например, $(w_k)$ с биортогональной системой
$(w'_k),$ то всякий элемент $w\in W$ разлагается в ряд $w=\sum_{k=1}^\infty \langle   w'_k, w\rangle w_k$
и для $x\in X$ и $T: X\to W$ получаем: $\langle   Tx, w'_k\rangle= \langle   x, T^*w'_k\rangle$ и, следовательно,
$Tx=\sum_{k=1}^\infty \langle   T^*w'_k, x\rangle w_k.$ Отметим, что в этом примере пространство $W$
сепарабельно, а пространство $X$ не обязано быть сепарабельным.

4)\,
Если пространство $X$ сепарабельно и обладает свойством ограниченной аппроксмации,
то оно дополняемо вкладывается в некоторое банахово пространство $X_0$ с базисом [3].       
Пусть $P: X_0\to X$ --- непрерывный проектоор и $((x'_k), (x_k))$ --- базис для $X_0$
(т.е., соответствующая биортогональная систеьа). Тогда для $x_0\in X_0$ мы имеем разложение
$x_0=\sum_{k=1}^\infty \langle  x'_k, x_0\rangle x_k$ и, если $x\in X$ и $x=Px_0,$ то
$x=\sum_{k=1}^\infty \langle  x'_k, x_0\rangle Px_k.$ Но $x\in X\subset X_0;$ следовательно,
$x=Px$ и  
$$
 x=\sum_{k=1}^\infty \langle x'_k, x\rangle x_k= P(x) =\sum_{k=1}^\infty \langle x'_k, x\rangle Px_k.
$$
Теперь, любой оператор $T\in L(X,W)$ имеет О-фрейм с разложением вида
$$
 Tx= \sum_{k=1}^\infty \langle x'_k|_X, x\rangle TPx_k.
$$
Аналогичный пример можно рассмотреть и для случая, когда $W$ сепарабельно
и обладает свойством ограниченной аппроксимации.
\smallskip

{\bf Лемма 1.1.}\, {\it
Пусть $X, V, Z, W$ --- банаховы пространства,
$T\in L(X,W),$ $A\in L(W,V),  B\in L(Z,X).$ Если оператор $T$
имеет О-фрейм, то и оператор $ATB: Z\to V$ имеет О-фрейм.
}
\smallskip

{\it Доказательство}.\,
Пусть $((x'_k), (w_k))$ --- О-фрейм для $T.$ Зафиксируем $z\in Z,$ и пусть
$x=Bz.$ Тогда $Tx=\sum \langle x'_k, x\rangle w_k,$ $ATx=\sum \langle x'_k, x\rangle Aw_k$ и
$ATBz=ATx=\sum \langle x'_k, Bz\rangle Aw_k= \sum \langle B^* x'_k,z\rangle Aw_k.$
Таким образом, $((B^*x'_k), (Aw_k))$ есть О-фрейм для оператора $ATB.$
\smallskip

{\bf Следствие 1.1.}\, {\it
Если оператор $T\in L(X,W)$ факторизуется через банахово пространство 
с базисом Шаудера, то он имеет О-фрейм.
}
\smallskip

{\it Доказательство}.\,
Применяем лемму 1.1 и примеры 2--3.
\smallskip

{\bf Следствие 1.2.}\, {\it
Если оператор $T\in L(X,W)$ факторизуется через банахово пространство 
со свойством ограниченной аппроксимации, то он имеет О-фрейм.
}
\smallskip

И еще одно свойство О-фреймов:
\smallskip

{\bf Предложение 1.1.}\,
Пусть $\mathcal F:=((x'_k), (w_k))$  есть О-фрейм для $T\in L(X.W).$
Тогда дуальная система $\mathcal F^d:=((w_k), (x'_k))$ есть $w^*$-слабый О-фрейм для $T^*,$
то есть 
$$
 T^*w'= w^*\text{-}\lim_N \sum_{k=1}^N \langle w', w_k\rangle x'_k,\ \ w'\in W^*.
$$
\smallskip

{\it Доказательство}.\,
Для $w'\in W^*$ и $x\in X$ имеем:
$$
 \langle Tx,w'\rangle = \langle \sum_{k=1}^\infty \langle x'_k, x\rangle w_k, w'\rangle =\langle \sum_{k=1}^\infty \langle w', w_k\rangle x'_k, x\rangle,
$$
откуда
$T^*w' =  w^*\text{-}\lim_N \sum_{k=1}^N  x'_k \langle w', w_k\rangle $
(предел в топологии $\sigma(X^*,X)$). 
\smallskip

Пока мы рассмотрели лишь простейшие свойства О-фреймов.
Намного более серьезные результаты появятся ниже.

В теории базисов и фреймов в банаховых пространствах часто возникают термины
"натягивающий", "ограниченно полный" (см., например, [1], [2]).     
Мы вводим соответствующие понятия для операторных фреймов.
\smallskip

{\bf Определение 1.2.}\,
О-фрейм $((x'_k), (w_k))$ для $T$ называется ограниченно полным, если
для любого $x''\in X^{**}$ ряд $\sum_{k=1}^\infty \langle x'', x'_k\rangle w_k$  
сходится в пространстве $W.$
\smallskip

{\bf Определение 1.3.}\,
О-фрейм $((x'_k), (w_k))$ для $T$ называется натягивающим, если
для любого $w'\in W^{*}$ 
норма  $||\sum_{k=n+1}^\infty  x'_k \langle w', w_k\rangle||\to 0$   при $n\to \infty.$

\smallskip

{\bf Предложение 1.2.}\,
Пусть $\mathcal F:=((x'_k), (w_k))$  есть О-фрейм для $T\in L(X.W).$
Дуальная система $\mathcal F^d:=((w_k), (x'_k))$ есть  О-фрейм для $T^*$
тогда и только тогда, когда
О-фрейм $\mathcal F$ натягивающий.
\smallskip

{\it Доказательство}.\,
Если О-фрейм $\mathcal F$ является натягивающим, то, во-первых (предложение 1.1),
для $w'\in W^*$ 
$$
 T^*w'= w^*\text{-}\lim_N \sum_{k=1}^N \langle w', w_k\rangle x'_k
 $$
 и, во-вторых,
$$||\sum_{k=n+1}^\infty  x'_k \langle w', w_k\rangle||\to 0, \ \ \text{ при\  }\ n\to\infty.$$
Следовательно, ряд $\sum_{k=1}^\infty \langle w',w_k\rangle x'_k$ сходится в $X^*,$ причем к $T^*w'.$

Обратно, если ряд $\sum_{k=1}^\infty \langle w',w_k\rangle x'_k$ сходится к $x_0$ в $X^*,$
то, по тому же предложению 1.1, он сходится к $T^*w'.$ 
Следовательно, остаток этого ряда стремится к нулю в $X^*,$ то есть,
О-фрейм $\mathcal F$ натягивающий.
\smallskip

{\bf Предложение 1.3.}\,
Пусть $\mathcal F:=((x'_k), (w_k))$  есть О-фрейм для $T\in L(X.W).$
Следующие утверждения равносильны:

$1)$\, О-фрейм $\mathcal F$ ограниченно полный;

$2)$\, для каждого $x''\in X^{**}$ из ограниченности последовательности частичных сумм
$(\sum_{k=1}^N \langle x'', x'_k\rangle w_k)_{N=1}^\infty$ вытекает сходимость ряда
$\sum_{k=1}^\infty \langle x'', x'_k\rangle w_k$ в пространстве $W.$
\smallskip

{\it Доказательство}.\,
Импликация $1)\Rightarrow 2)$ очевидна.

$2)\Rightarrow 1).$\,
Пусть $x''\in X^{**}.$
Надо показать, что ряд $\sum_{k=1}^\infty \langle x'', x'_k\rangle w_k$ сходится в $W.$

Для любых $w'\in W^*$ и $x\in X$ имеем:
$$
 \lim_{N\to\infty} \langle \sum_{k=1}^N \langle w', w_k\rangle x'_k, x\rangle =
 \lim_{N\to\infty} \langle w', \sum_{k=1}^N  w_k\langle  x'_k, x\rangle\rangle= \langle w', Tx\rangle.
$$
В частности, последовательность $(\sum_{k=1}^N \langle w', w_k\rangle x'_k)_{N=1}^\infty$
ограничена в $X^*.$ 
Поэтому, для $w'\in W^*$ и $x''\in X^{**}$ числовая последовательность
$$(\sum_{k=1}^N \langle w', w_k\rangle  \langle x'', x'_k\rangle)_{N=1}^\infty$$
ограничена, то есть последовательность
$(\sum_{k=1}^N  w_k  \langle x'', x'_k\rangle)_{N=1}^\infty$ является
$\sigma(W,W^*)$-ограниченной. Следовательно,
она сильно ограничена в $W.$ По условию 2), ряд
$\sum_{k=1}^\infty \langle x'', x'_k\rangle w_k$ сходится в пространстве $W,$
т.е. наш О-фрейм $\mathcal F$ --- ограниченно полный.
\smallskip

Применяя оба предложения, получаем один из основных фактов данного раздела:

\smallskip

{\bf Теорема 1.1.}\,
Пусть $\mathcal F:=((x'_k), (w_k))$  есть О-фрейм для $T\in L(X.W).$
Если этот  О-фрейм $\mathcal F$ ограниченно полный и натягивающий,
то оператор $T$ слабо компактен.
\smallskip

{\it Доказательство}.\,
Так как $\mathcal F$ натягивающий, то по предложению 1.2
$\mathcal F^d$ есть О-фрейм для $T^*,$ т.е.
для любого $w'\in W^*$
$$
 T^*w' =\sum_{k=1}^\infty \langle w', w_k\rangle x'_k.
$$
Поэтому, для $x''\in X^{**}$ 
$$
 \langle T^*w', x''\rangle = \sum_{k=1}^\infty \langle w', w_k\rangle \langle x'', x'_k\rangle= \lim_N \langle w', \sum_{k=1}^N \langle x'', x'_k\rangle w_k\rangle,
$$
или
$$
 \langle w', T^{**}x''\rangle =  \lim_N \langle w', \sum_{k=1}^N \langle x'', x'_k\rangle w_k\rangle.
$$
Так как $\mathcal F$ есть ограниченно полный О-фрейм, то
ряд $\sum_{k=1}^\infty \langle x'', x'_k\rangle w_k$ сходится в $W,$ пусть, например, к элементу $w_0\in W.$
Следовательно, 
$$
 \langle w', T^{**}x''\rangle =   \langle w', \sum_{k=1}^\infty \langle x'', x'_k\rangle w_k\rangle. \ \ w'\in W^*.
$$
Или: $\langle T^{**}x'', w'\rangle = \langle w', w_0\rangle$ для всякого $w'\in W^*.$
Таким образом, $T^{**}x''=w_0$ и, в частности, $T^{**}x''\in W,$
то есть, если $x''\in X^{**},$ то $T^{**}x''\in W$ и, следовательно, оператор
$T$ слабо компактен.
\smallskip


Следующая теорема --- первая из теорем, которые связывают введенное нами 
понятие О-фрейма и понятие свойства BAP
(ограниченной аппроксимации) для операторов (подробней о BAP см. ниже).
\smallskip
 
{\bf Теорема 1.2.}\, 
{\it
Пусть $T\in L(X,W).$ Следующие утверждения равносильны:

$1)$\,
$T$ имеет О-фрейм;

$2)$\,
оператор $T$ факторизуется через банахово пространство с базисом Ша\-удера;

$3)$\,
оператор $T$ факторизуется через банахово пространство последова\-тель\-ностей 
с базисом Шаудера;
}
\smallskip

{\it Доказательство}.\,
Сначала мы покажем, что из 1) вытекает 3).
Пусть оператор $T$ имеет О-фрейм $\mathcal F:=((x'_k), (w_k)).$ Мы можем (и будем)
считать, что $w_k\neq0$ для каждого $k.$ 
Так как ряд $\sum_{k=1}^\infty \langle x'_k, x\rangle w_k$ сходится (к $Tx)$ для каждого $x\in X,$
то нормы операторов $\sum_{k=1}^N x'_k\otimes w_k$ (равномерно) ограничены (например,
константой $K>0).$

Рассмотрим пространство числовых последовательностей 
$$
 t:=\{a=(a_k)_{k=1}^\infty:\  \text{ряд}\ \ \sum_{k=1}^\infty a_kw_k\ \  \text{сходится\ в}\ W\},   
$$
и пусть $e_k$ --- $k$-й единичный вектор в $t$
(т.е., $(e_k)_s=0$ при $k\neq s$ и $(e_k)_k=1).$

Положим
$$
  |||a|||_t:= \sup_N ||\sum_{k=1}^N a_kw_k||\ \, (\ge \lim_N ||\sum_{k=1}^N a_kw_k||).
$$
Так как для финитной последовательности $a=(a_1, a_2,\dots,  a_{N+s}, 0, 0,\dots)$
$$
 |||\sum_{k=1}^N a_ke_k||| \le |||\sum_{k=1}^{N+s} a_ke_k|||
$$
и линейная оболочка векторов $(e_k)_{k=1}^\infty$ плотна в $t,$ то $(e_k)$ --- монотонный базис
(см. [2]) в                                                
банаховом пространстве $t.$ 
При этом, если $j: t\to W$ --- естественное отображение
$a\mapsto \sum_{k=1}^\infty a_kw_k,$ то $||j||\le1$ 
(по определению нормы $|||\cdot|||$ в $t).$
Положим $Ax:= (\langle x'_k,x\rangle)_{k=1}^\infty;$ так как ряд
$\sum_{k=1}^\infty \langle x'_k, x\rangle w_k$ сходится,  то $Ax\in t.$
Далее,
$$
 |||Ax|||_t=\sup_N ||\sum_{k=1}^N \langle x'_k, x\rangle w_k||\le K\, ||x||,\ \ \forall\, x\in X.
$$
Таким образом, линейное отображение $A$ непрерывно из $X$ в $t,$
и $T=jA: X\to t\to W.$

3) влечет 2) --- очевидно, а утверждение 1) вытекает из утверждения 2) по следствию 1.1.
\smallskip

Наряду с понятием О-фрейма, мы можем рассмотреть (и кратко рассмотрим)
понятие безусловного О-фрейма.
\smallskip

{\bf Определение 1.4.}\,
Пусть $T\in L(X,W),$ $(x'_k)_{k=1}^\infty\subset X^*, (w_k)_{k=1}^\infty\subset W.$  
Мы говорим, что $\mathcal F:=((x'_k)_{k=1}^\infty, (w_k)_{k=1}^\infty)$
есть БО-фрейм для $T,$ если для каждого $x\in X$
ряд $\sum_{k=1}^\infty \langle x'_k, x\rangle w_k$ сходится безусловно в $W$ и
$$
  Tx= \sum_{k=1}^\infty \langle x'_k, x\rangle w_k, \ \ x\in X.
$$
Если существует БО-фрейм для $T,$ то мы говорим, что $T$ имеет БО-фрейм.
\smallskip

{\bf Теорема 1.3.}\, 
{\it
Пусть $T\in L(X,W).$ Следующие утверждения равносильны:

$1)$\,
$T$ имеет БО-фрейм;

$2)$\,
оператор $T$ факторизуется через банахово пространство с безусловным базисом;

$3)$\,
оператор $T$ факторизуется через банахово пространство последовате\-ль\-ностей 
с безусловным базисом.
}
\smallskip

{\it Доказательство}.\,
Сначала мы покажем, что из 1) вытекает 3).
Пусть оператор $T$ имеет БО-фрейм $\mathcal F:=((x'_k), (w_k)).$ Мы можем (и будем)
считать, что $w_k\neq0$ для каждого $k.$ 
Так как ряд $\sum_{k=1}^\infty \langle x'_k, x\rangle w_k$ безусловно сходится (к $Tx)$ для каждого $x\in X,$
то  билинейная форма $B(\cdot,\cdot): X\times l_\infty \to X,$ определенная
для $x\in X$ и $a\in l_\infty$ формулой
$B(x,a):= \sum_{k=1}^\infty a_k\langle x'_k, x\rangle w_k,$ ограничена.
 
Рассмотрим пространство числовых последовательностей 
$$
u:=\{a=(a_k)_{k=1}^\infty:\  \text{ряд}\ \ \sum_{k=1}^\infty a_kw_k\ \  \text{безусловно\ сходится\ в}\ W\},   
$$
и пусть $e_k$ --- $k$-й единичный вектор в $t$
(т.е., $(e_k)_s=0$ при $k\neq s$ и $(e_k)_k=1).$

Положим
$$
  |||a|||_u:= \sup_{||(b_k)||_{l_\infty}\le1} ||\sum_{k=1}^\infty b_ka_kw_k||.
$$
Тогда $(e_k)$ --- безусловный базис в                                                
банаховом пространстве $u.$ 
При этом, если $j: u\to W$ --- естественное отображение
$a\mapsto \sum_{k=1}^\infty a_kw_k,$ то $||j||\le1$ 
(по определению нормы $|||\cdot|||$ в $u).$
Положим $Ax:= (\langle x'_k,x\rangle)_{k=1}^\infty;$ так как ряд
$\sum_{k=1}^\infty \langle x'_k, x\rangle w_k$ безусловно сходится,  то $Ax\in t.$
Далее,
$$
 |||Ax|||_t=\sup_{||(b_k)||_{l_\infty}\le1} ||\sum_{k=1}^\infty b_k \langle x'_k, x\rangle w_k||\le ||B||\, ||x||
     \ \ \forall\, x\in X.
$$
Таким образом, линейное отображение $A$ непрерывно из $X$ в $u,$
и $T=jA: X\to u\to W.$

То, что 3) влечет 2) и 2) влечет 1) --- очевидно.
\smallskip

В качестве одного из следствий последнего факта, приведем ответ на вопрос Марии Скопиной
(примерно трех-летней давности; тогда я ответил на вопрос о существовании в $L_1$
безусловной системы представления, используя совершенно другие соображения):

{\bf Следствие 1.3.}\,
Пространство $L_1(0,1)$ не имеет счетной безусловной системы представления,
то есть, тождественный оператор в этом пространстве не обладает безусловным фреймом.

Для доказательства достаточно вспомнить теорему А. Пелчинского (см. [2, p.24; 1.d.1]):
пространство $L_1(0,1)$ не является подпространством банахова пространства с безусловным базисом.

\smallskip
Наша следующая цель --- связать понятие О-фрейма со свойством ограниченной 
аппроксимации операторов. 
\smallskip

\vskip 0.3cm

\centerline{\bf \S2. Об операторах со свойством ограниченной аппроксимации}

 \vskip 0.23cm

{\bf Определение 2.1.}\,
Пусть $T\in L(X,W),$  $C\ge1.$ Мы говорим, что $T$ имеет свойство C-BAP (свойство C-ограниченной аппроксимации),
если для всякого компактного подмножества $K$ из $X,$  для каждого $\varepsilon>0$ существует такой конечномерный оператор
 $R: X\to W,$ что  $||R||\le C\, ||T||$ и $\sup_{x\in K} ||Rx-Tx||\le \varepsilon.$ 
Оператор $T$ имеет (свойство) BAP, если он имеет свойство C-BAP для некоторой
постоянной $C\in[1,\infty).$


\medskip 

{\bf Лемма 2.1.}\,
$T$ имеет C-BAP тогда и только тогда, когда для любого конечного семейства
 $(x_k)_{k=1}^N\subset X,$ для всякого $\varepsilon>0$ существует такой конечно\-мерный оператор
 $R: X\to W,$ что
$||R||\le C||T||$ и $\sup_{1\le k\le N} ||Rx_k-Tx_k||\le \varepsilon.$

{\it Доказательство}.\,
Мы можем (и будем) считать, что $||T||=1.$
Зафиксируем компактное подмножество $K\subset X$ и $\varepsilon>0.$ Пусть $\varepsilon_0:= \varepsilon/(2+C),$ $(x_k)_{k=1}^M$ ---
 $\varepsilon_0$-сеть для $K$ в $X,$ $R\in X^*\otimes W,$ $||R||\le C$ и $\sup_{1\le k\le M} ||Rx_k-Tx_k||\le \varepsilon_0.$
Возьмем произвольный $x\in K,$ и пусть  $x_k$ таков, что $||x-x_k||\le \varepsilon_0.$  Тогда
$||Tx-Rx||\le ||x-x_k||+ ||Tx_k-Rx_k|| +||Rx_k-Rx||\le \varepsilon_0+\varepsilon_0+C\varepsilon_0=\varepsilon.$
\smallskip

{\bf Лемма 2.2.}\,
Пусть $X, W$ --- банаховы пространства, причем пространство  $X$ сепарабельно, и $T\in L(X, W).$
$T$ имеет C-BAP тогда и только тогда, когда существует последовательность $(Q_l)_{l=1}^\infty$ 
конечномерных операторов из $X$ в $W,$ для которой

$1)$\, при любом $x\in X$ ряд $\sum_{l=1}^\infty Q_lx$ сходится и
$$
 Tx=\sum_{l=1}^\infty Q_lx,\ x\in X;
$$

$2)$\, $\sup_N ||\sum_{l=1}^N Q_l||\le C\, ||T||.$
\smallskip

{\it Доказательство}.\,
Так как $X$ сепарабельно, то существует последователь\-ность $(x_k)_1^\infty,$ которая плотна
в замкнутом единичном шаре $\bar B_1(0)$ пространства $X.$
Предположим, как и выше, что $||T||=1$ и $T$ имеет C-BAP, то есть
для каждого конечного множества $F\subset X,$ для любого $\varepsilon>0$ существует такой конечномерный 
оператор $R:X\to W,$ что $||R||\le C$ и $\sup_{f\in F} ||Rf-Tf||\le \varepsilon.$
Для  каждого $N,$ пусть $R_N$ есть конечнмерный оператор из $X$ в $W$ со следующими свойствами:

(i)\, $||R_N||\le C$ и

(ii)\, $\sup_{1\le n\le N} ||R_Nx_n-Tx_n||\le 1/2^{N+1}.$

Если $n\in \Bbb N,$ то для всякого $N\ge n$ имеем:
$$
 (iii)\ \ ||R_Nx_n-Tx_n||\le \frac1{2^{N+1}}
$$
и, следовательно, для фиксированного $x_n$
$$
 R_Nx_n \rightarrow Tx_n
$$
при $N$ стремящемся к $\infty.$

Теперь, зафиксируем $\varepsilon>0$ и пусть $\delta>0$ таково, что $C\delta+\delta<\varepsilon.$
Для $x\in \bar B_1(0),$  возьмем такой $x_n,$ что  $||x_n-x||<\delta.$
Тогда найдется такой номер $N_0,$ что при $N\ge N_0$
$$
 ||R_Nx-Tx||\le ||R_N||\, ||x_n-x|| +||R_Nx_n-Tx_n||+ ||Tx_n-Tx||\le C\delta +||T||\delta<\varepsilon.
$$
Таким образом, если $x\in X,$ то $R_Nx \to Tx$ при $N\to +\infty.$

Для завершения доказательства в части "только тогда"  мы применим лемму:
\smallskip

{\bf Лемма 2.3.}\,
Пусть $X, W$ --- произвольные банаховы пространства,  $C\ge1$ и $T\in L(X, W).$
Предположим, что

(*)\, существует такая последовательность $(S_N)_{N=1}^\infty$ конечномерных операторов
из $X$ в $W,$ что
если  $x\in X,$ то $S_Nx \to Tx$ при $N\to +\infty$ и $||S_N||\le C||T||$ для каждого $N.$

Тогда существует такая последовательность
 $(Q_l)_{l=1}^\infty$ конечномерных операторов из $X$ в $W,$ что

$1)$\, для любого $x\in X$ ряд $\sum_{l=1}^\infty Q_lx$ сходится и
$$
 Tx=\sum_{l=1}^\infty Q_lx,\ x\in X;
$$

$2)$\, $\sup_N ||\sum_{l=1}^N Q_l||\le C\, ||T||.$
\smallskip

{\it Доказательство}.\,
Мы снова предполагаем, что  $||T||=1.$
Положим  $Q_1:=S_1, Q_l:=S_l-S_{l-1}$ для $l>1,$ так что
$$
  S_N= S_1+(S_2-S_1)+\dots +(S_{N-1}-S_{N-2})+ (S_N-S_{N-1})=
  Q_1+Q_2+\dots+Q_N.
$$
Отсюда следует, что

$$(1)\ Tx=\sum_{l=1}^\infty Q_lx\ \ \forall\, x\in X$$
и
$$(2)\ \sup_{N}||\sum_{l=1}^NQ_l||=\sup_N ||S_N||\le C.$$
\smallskip

Часть "тогда, когда" доказательства леммы 2.2 вытекает из следующего утверждения.
\smallskip

{\bf Лемма 2.4.}\,
Пусть $X, W$ --- произвольные банаховы пространства, $C\ge1$ и $T\in L(X,W).$
Если существует последовательность $(R_N)_1^\infty$
конечномерных операторов из $X$ в $W,$  сходящаяся поточечно к $T$
и такая. что $||R_N||\le C||T||$ для всех  $N,$
то $T$ имеет C-BAP.
\smallskip

{\it Доказательство}.\,
Действительно, пусть $R_Nx\to Tx$ для каждого $x\in X$ (и $||R_N||\le C||T||$). Зафиксируем  $\varepsilon>0$
и компактное подмножество $K\subset X.$ Положим $\varepsilon_0:=\varepsilon (||T||+1+C||T||)^{-1}.$
Выберем конечную $\varepsilon_0$-сеть $F\subset X$ для $K$ и рассмотрим такой оператор  $R_{N_0},$ что
$\sup_{f\in F} ||R_{N_0}f- Tf||\le \varepsilon_0.$ Тогда, для любого $x\in K$ найдется такой элемент
 $f_0\in F,$ что $||f_0-x||\le\varepsilon_0,$ и мы получаем:
$$
  ||Tx-R_{N_0}x||\le ||T||\, \varepsilon_0+\varepsilon_0+||R_{N_0}||\varepsilon_0\le
  \varepsilon_0 (||T||+1+C||T||)=\varepsilon.
$$
\smallskip

{\bf Следствие 2.1.}\,
Если $X$ сепарабельно и $T\in L(X,W),$ то $T$
имеет C-BAP тогда и только тогда, когда существует последовательность $(R_N)_1^\infty$
конечномерных операторов из  $X$ в $W,$ которая поточечно сходится к оператору $T$
и для которой $||R_N||\le C||T||$ для всех $N.$
\smallskip

{\bf Следствие 2.2.}\,
Если $X$ сепарабельно и $T\in L(X,W),$ то $T$
имеет BAP тогда и только тогда, когда существует последовательность 
конечномерных операторов из  $X$ в $W,$ сходящаяся к оператору $T$ поточечно.

\bigskip

\centerline{\bf \S3. Свойство ограниченной аппроксимации и О-фреймы.}

\bigskip

В этом разделе мы покажем, в частности, что каждый оператор в сепарабе\-льных
банаховых пространствах, обладающий свойством ограниченной апро\-к\-си\-мации, имеет О-фрейм.
 
\smallskip

Мы будем использовать результаты предыдущего параграфа, но заменим некоторые обозначения,
поскольку основная идея доказательства центральной теоремы 3.1 восходит к А. Пелчинскому и
нам представляется, что наилучшим вариантом является вариант с использованием 
обозначений из его важной рабо\-ты [3].                         

 Пусть $X, W$ --- произвольные банаховы пространства,
$T\in L(X,W)$ и $T$ обладает свойством (*) из Леммы 2.3.
 Рассмотрим последовательность $(Q_l)_{l=1}^\infty,$ полученную в утверждении Леммы 2.3,
 и положим
 $A_p:=Q_p$\, $(p=1,2,\dots)$ и $K:=C (\ge1),$ предполагая, что $||T||=1.$
 Сейчас мы --- в обозначениях (частично) статьи [3].          
\smallskip

{\bf Теорема 3.1.}\,
Если $T:X\to W$ имеет свойство (*), то есть, если
существует такая последовательность $(S_N)_{N=1}^\infty$ конечномерных операторов
из $X$ в $W,$ что
для любого   $x\in X$ \  $S_Nx \to Tx$ при $N\to +\infty$ и $||S_N||\le C||T||$ для каждого $N,$
то оператор $T$ имеет О-фрейм.
\smallskip

{\it Доказательство}.\,
В вышеприведенных обозначениях (предполагая, что $||T||=1$), имеем:
$$
 Tx=\sum_{p=1}^\infty A_px,\ \forall\, x\in X;\  \ A_p\in X^*\otimes W,\  \sup_{n\in\Bbb N} ||\sum_{p=1}^n A_p||\le K
$$
(заметим, что для каждого $n$ \, $||A_n||\le ||\sum_{p=1}^n A_p- \sum_{p=1}^{n-1} A_p||\le 2K$).
Пусть $E_p=A_p(X)\subset W,$ \, $m_p:=\operatorname{dim} E_p$ для $p\ge1$ и $m_0=0.$
Мы будем действовать как в [3].                     

По лемме Ауэрбаха, существуют такие одномерные операторы
$B_j^{(p)}: E_p\to E_p,$ что
 $||B_j^{(p)}||=1$ для $j=1,2,\dots, m_p,$
и
$$
 \sum_{j=1}^{m_p} B_j^{(p)}(e)=e,\ e\in E_p.
$$
Положим  $C_i^{(p)}:= \dfrac1{m_p}\, B_j^{(p)}$ для $i=rm_p+j$ (где $r=0,1,\dots, m_p-1; j=1,2,\dots, m_p).$
Тогда, для $e\in E_p$
$$
  \sum_{i=1}^{m_p^2} C_i^{(p)} (e)=m_p \cdot \sum_{j=1}^{m_p} \frac1{m_p} B_j^{(p)}e=
  \sum_{j=1}^{m_p}  B_j^{(p)}e=e.
$$
 Кроме того, для любого $q\ge1, q\le m_p^2$ и некоторых $l< m_p$ и $k\le m_p$ имеем:
 $$
   ||\sum_{i=1}^q C_i^{(p)}||= ||\sum_{i=1}^{lm_p} C_i^{(p)} +\sum_{lm_p+1}^{lm_p+k} C_i^{(p)}|| \le
       l\cdot \frac1{m_p} \, ||\sum_{j=1}^{m_p} B_j^{(p)}|| + \frac1{m_p}\cdot ||\sum_{j=1}^{k} B_j^{(p)}|| \le 1+1=2.
 $$

Теперь, пусть
 $$\widetilde A_s:= C_i^{(p)}A_p$$
  для $p\in \Bbb N,$ $i=1,2,\dots, m_p^2$ и $s=m_0^2+m_1^2+\dots +m_{p-1}^2+i.$
Одномерный оператор $\widetilde A_s$ отображает $X$ в $E_p\subset W$ следующим образом:
  $$
   \widetilde A_s: X\stackrel{A_p}{\rightarrow} E_p=A_p(X)\stackrel{C_i^{(p)}}{\rightarrow} E_p (\subset W).
$$
 Так как для каждого $n\in\Bbb N,$ для некоторых $k$ и $r\le m_k^2$

 $$
  \sum_{s=1}^n \widetilde A_s=\sum_{p=1}^{k-1} \sum_{i=1}^{m_p^2} C_i^{(p)}A_p + \sum_{i=1}^{r} C_i^{(k)}A_k,
 $$
то мы получаем, что
 $$
  || \sum_{s=1}^n \widetilde A_s||\le ||\sum_{p=1}^{k-1} A_p|| + ||\sum_{i=1}^{r} C_i^{(k)}A_k||\le K+2||A_k||\le 5K.
 $$
 
 Поскольку для каждого $x,\, x\in X,$\, $A_kx\to 0$ при $k\to\infty,$ то:
 $$
  \lim_{n\to\infty}  \sum_{s=1}^n \widetilde A_s x= \lim_{N\to\infty} \sum_{p=1}^{N} A_px =Tx.
 $$
Таким образом, $\sum_{s=1}^\infty \widetilde A_s(x)= Tx,$ где
$\widetilde A_s\in X^*\otimes W$ --- одномерный оператор для любого $s,$
откуда и следует, что оператор $T$ имеет О-фрейм.

\medskip

Из доказанной теоремы и результатов предыдущих параграфов получаем важнейшие следствия.

\smallskip

{\bf Следствие 3.1.}\,
Для любых банаховых пространств $X$ и $W$ и для любого оператора $T\in L(X,W)$
равносильны утверждения:

$(1)$\, $T$ имеет О-фрейм;

$(2)$\, оператор $T$ факторизуется через банахово пространство 
последователь\-но\-стей с базисом;

$(3)$\, существует последовательность 
        конечномерных операторов из  $X$ в $W,$ сходящаяся к оператору $T$ поточечно.

\smallskip

{\bf Следствие 3.2.}\,
Пусть $X$ --- сепарабельное банахово пространство, $W$ --- произвольное банахово
пространство и $T\in L(X,W).$ Следующие утверждения эквивалентны:

$(1)$\, $T$ имеет О-фрейм;

$(2)$\, $T$ обладает свойством ограниченной аппроксимации;

$(3)$\, оператор $T$ факторизуется через банахово пространство 
         с базисом Шау\-де\-ра.

\smallskip



\end{document}